\documentclass[12pt]{amsart}
    \title[Compactly Generated Co-$t$-structures]{A Note on Compactly Generated Co-$t$-structures}
    \author{David Pauksztello}
    \date{12th October 2010}
    
\usepackage{amssymb}
\usepackage{latexsym}
\usepackage{amsmath}
\usepackage{mathrsfs}
\usepackage[all]{xy}

\newtheorem{dfn}{Definition}
\newtheorem{thrm}[dfn]{Theorem}

\newtheorem{lem}[dfn]{Lemma}
\newtheorem{prp}[dfn]{Proposition}
\newtheorem{corr}[dfn]{Corollary}

\newtheorem{rmk}[dfn]{Remark}
\newtheorem{rmks}[dfn]{Remarks}

\newtheorem{set}[dfn]{Setup}
\newtheorem{example}[dfn]{Example}
\newtheorem{examples}[dfn]{Examples}
\newtheorem{observe}[dfn]{Observation}
\newtheorem{definitions}[dfn]{Definitions}
\newtheorem{conject}[dfn]{Conjecture}
\newenvironment{prf}{\noindent\textbf{Proof. }}{\hfill $\Box$\ }
\newenvironment{prfoflemma}{\noindent\textbf{Proof of lemma. }}{\hfill $\Box$\ }
\newenvironment{enumerateroman}{\begin{enumerate}}{\end{enumerate}}

\newcommand{\df}[1]{\begin{dfn}\emph{#1}}
\newcommand{\edf}{\end{dfn}}
\newcommand{\thm}{\begin{thrm}}
\newcommand{\ethm}{\end{thrm}}
\newcommand{\lemma}{\begin{lem}}
\newcommand{\elemma}{\end{lem}}
\newcommand{\prop}{\begin{prp}}
\newcommand{\eprop}{\end{prp}}
\newcommand{\rem}[1]{\begin{rmk}\emph{#1}}
\newcommand{\erem}{\end{rmk}}
\newcommand{\rems}[1]{\begin{rmks}\emph{#1}}
\newcommand{\erems}{\end{rmks}}
\newcommand{\pf}{\begin{prf}}
\newcommand{\epf}{\end{prf}}
\newcommand{\cor}{\begin{corr}}
\newcommand{\ecor}{\end{corr}}
\newcommand{\setup}[1]{\begin{set}\emph{#1}}
\newcommand{\esetup}{\end{set}}
\newcommand{\eg}[1]{\begin{example}\emph{#1}}
\newcommand{\eeg}{\end{example}}
\newcommand{\egs}[1]{\begin{examples}\emph{#1}}
\newcommand{\eegs}{\end{examples}}
\newcommand{\obs}[1]{\begin{observe}\emph{#1}}
\newcommand{\eobs}{\end{observe}}
\newcommand{\dfs}[1]{\begin{definitions}\emph{#1}}
\newcommand{\edfs}{\end{definitions}}
\newcommand{\conj}{\begin{conject}}
\newcommand{\econj}{\end{conject}}
\newcommand{\pflemma}{\begin{prfoflemma}}
\newcommand{\epflemma}{\end{prfoflemma}}

\renewcommand{\geq}{\geqslant}
\renewcommand{\leq}{\leqs}
\renewcommand{\phi}{\varphi}
\DeclareMathAlphabet{\mathpzc}{OT1}{pzc}{m}{it}
\newcommand{\D}{\mathscr{D}}

\newcommand{\T}{\mathscr{T}}
\newcommand{\X}{\mathscr{X}}
\newcommand{\Y}{\mathscr{Y}}
\renewcommand{\S}{\mathscr{S}}
\newcommand{\A}{\mathscr{A}}
\newcommand{\B}{\mathscr{B}}
\newcommand{\C}{\mathscr{C}}

\newcommand{\F}{\mathscr{F}}
\newcommand{\G}{\mathscr{G}}

\newcommand{\R}{\mathscr{R}}

\newcommand{\Z}{\mathbb{Z}}

\newcommand{\Hom}[3]{\mathsf{Hom}_{#1}(#2,#3)}

\newcommand{\rightiso}{\stackrel{\sim}{\longrightarrow}}

\newcommand{\hocolim}{\underrightarrow{\mathsf{holim}}\,}
\newcommand{\colim}{\underrightarrow{\mathsf{lim}}\,}
\newcommand{\limit}{\underleftarrow{\mathsf{lim}}\,}

\newcommand{\Endo}[1]{\mathsf{End}(#1)}

\newcommand{\Add}[1]{\mathsf{Add}(#1)}

\newcommand{\rightlabel}[1]{\stackrel{#1}{\longrightarrow}}

\newcommand{\trilabel}[4]{#1\stackrel{#4}{\longrightarrow} #2\longrightarrow #3\longrightarrow \Sigma #1}
\newcommand{\trilabels}[6]{#1\stackrel{#4}{\longrightarrow} #2\stackrel{#5}{\longrightarrow} #3\stackrel{#6}{\longrightarrow} \Sigma #1}
\newcommand{\backtri}[3]{\Sigma^{-1}#1\rightarrow #2\rightarrow #3\rightarrow #1}

\newcommand{\leqs}{\leqslant}

\entrymodifiers={+!!<0pt,\fontdimen22\textfont2>}

\begin{document}

\begin{abstract}
The idea of a co-$t$-structure is almost `dual' to that of a $t$-structure, but with some important differences. This note establishes co-$t$-structure analogues of Beligiannis and Reiten's corresponding results on compactly generated $t$-structures.
\end{abstract}
	
\address{Institut f\"ur Algebra, Zahlentheorie und Diskrete Mathematik, Fakult\"at f\"ur Mathematik und Physik, Leibniz Universit\"at Hannover, Welfengarten 1, 30167 Hannover, Germany.}
\email{pauk@math.uni-hannover.de}

     
    
\maketitle

The notion of a co-$t$-structure on a triangulated category was introduced independently by the author in \cite{Co-t-structures} and Bondarko in \cite{Bondarko}. In \cite{Bondarko} they are referred to as weight structures; in this note we continue the terminology of \cite{Co-t-structures}. In \cite{Bondarko} they are introduced in the context of understanding Grothendieck's weight filtration in Voevodsky's triangulated category of motives; see also \cite{Bondarko-coniveau}. In \cite{Achar_Treumann} co-$t$-structures are important ingredients in the proofs of purity and decomposition theorems for `staggered sheaves', and in \cite{Wildeshaus2} and \cite{Wildeshaus} they are studied in the context of Chow motives and Artin-Tate motives, respectively. In the representation theoretic setting, co-$t$-structures have recently been studied in connection with the Auslander-Buchweitz context in \cite{Mendoza}.

In view of their recent proliferation into different branches of mathematics, it is useful to establish basic results regarding the structure and behaviour of co-$t$-structures, and, in particular, their relation with, similarities to, and differences from $t$-structures. In this sense, the present note should be viewed as an extension of \cite{Co-t-structures} providing co-$t$-structure analogues of the corresponding results for $t$-structures in \cite{Beligiannis_Reiten}. 

Throughout this note, $\T$ will be a triangulated category with set indexed coproducts and and $\Sigma:\T\rightarrow\T$ will denote its suspension functor. We direct the reader to \cite{Holm_Jorgensen} for an introduction to triangulated categories. We first recall some definitions. An object $S$ of $\T$ is called {\it rigid} if $\Hom{\T}{S}{\Sigma^{i}S}=0$ for all $i>0$ (see \cite{Iyama}); $S$ is called {\it compact} if for any set indexed family of objects $\{X_{i}\}_{i\in I}$ of $\T$ one has the natural isomorphism $\Hom{\T}{S}{\coprod_{i\in I}X_{i}}\cong\coprod_{i\in I}\Hom{\T}{S}{X_{i}}$. Recall also that a {\it generating set} for $\T$ is a set of objects $\G$ such that if $\Hom{\T}{G}{X}=0$ for all $G\in\G$ then $X=0$.

In \cite[Theorem 1.3]{Miyachi}, if one takes a compact rigid object $S$ of $\T$ such that $\{\Sigma^{i}S\,|\,i\in\Z\}$ is a generating set, then one obtains a canonical $t$-structure on $\T$ given by
\begin{eqnarray}
\X & = & \{X\in\T \ | \ \Hom{\T}{S}{\Sigma^{i}X}=0 \textnormal{ for } i>0\}, \label{BR_t-structure} \\
\Y & = & \{X\in\T\ | \ \Hom{\T}{S}{\Sigma^{i}X}=0 \textnormal{ for } i<0\}. \nonumber
\end{eqnarray}

In \cite[Theorem III.2.3]{Beligiannis_Reiten}, it is established that any compact object $S$ of $\T$ induces a canonical $t$-structure on $\T$ given by
\begin{eqnarray}
\X & = & \Sigma^{-1}({}^{\perp}\Y), \label{BR_t-structure2} \\
\Y & = & \{X\in\T\ | \ \Hom{\T}{S}{\Sigma^{i}X}=0 \textnormal{ for } i<0\}. \nonumber
\end{eqnarray}
Since the `torsion-free class' $\Y$ of a $t$-structure is always closed under non-positive suspensions, this description of $\Y$ is natural. However, an obvious question is whether the `torsion class' $\X$ also has such a nice description. In \cite[Proposition III.2.8]{Beligiannis_Reiten}, this is shown to be the case if and only if the hypotheses of \cite[Theorem 1.3]{Miyachi} hold. A natural question is thus: what happens if one specifies the `torsion class' $\X=\{X\in\T \ | \ \Hom{\T}{S}{\Sigma^{i}X}=0 \textnormal{ for } i>0\}$ and sets $\Y=\Sigma(\X^{\perp})$? Unfortunately, it seems that this is not possible in general for $t$-structures but the situation arises naturally in the setting of co-$t$-structures.

Another motivation is that the theory of co-$t$-structures seems to be richer when there exist adjacent $t$-structures (see \cite[Sections 4.4 and 4.5]{Bondarko}). However, as remarked in \cite[Remark 4.5.3]{Bondarko}, the question of existence of an adjacent co-$t$-structure when a triangulated category is endowed with a $t$-structure is difficult in general. The main result of this note provides a case where such adjacent (co)-$t$-structures exist.

We now recall the definition of a co-$t$-structure:

\newtheorem{co-t-struct}[dfn]{Definition}
\begin{co-t-struct}[\cite{Co-t-structures}, Definition 2.4]
\textnormal{Let $\T$ be a triangulated category with set indexed coproducts. A pair of full subcategories of $\T$, $(\A,\B)$, is called a {\it co-$t$-structure} on $\T$ if it satisfies the following properties:
\begin{itemize}
\item[(0)] $\A$ and $\B$ are closed under direct summands;
\item[(1)] $\Sigma^{-1}\A\subseteq \A$ and $\Sigma\B\subseteq \B$;
\item[(2)] $\Hom{\T}{\Sigma^{-1}\A}{\B}=0$;
\item[(3)] For any object $X$ of $\T$ there exists a distinguished triangle $\backtri{A}{X}{B}$ with $A\in\A$ and $B\in\B$.
\end{itemize}}
\label{co-t-structure}
\end{co-t-struct}

In \cite[Definition 3.2]{Co-t-structures} an object $S$ of $\T$ is called a {\it simply connected corigid object} of $\T$ if
\begin{itemize}
\item[(1)] $S$ is {\it corigid}, that is, $\Hom{}{\Sigma^{i}S}{S}= 0$ for $i>0$;
\item[(2)] $\Hom{\T}{S}{\Sigma S}=0$;
\item[(3)] $\Endo{S}$ is a division ring.
\end{itemize}
We shall refer to $S$ as a {\it connected corigid object} of $\T$ if the first two conditions hold.

Motivated by the observation that the case studied in \cite{Beligiannis_Reiten} and \cite{Miyachi} represents a `chain situation', and that, often, a `cochain situation' is more natural (see \cite{Co-t-structures} for more remarks on this), the author showed in \cite[Theorem 3.2]{Co-t-structures} that if $S$ is a simply connected corigid object of $\T$ and $\{\Sigma^{i}S\,|\,i\in\Z\}$ is a generating set for $\T$, then one has a canonical co-$t$-structure on $\T$ given by 
\begin{eqnarray*}
\A & = & \{X\in\T \ | \ \Hom{\T}{S}{\Sigma^{i}X}=0 \textnormal{ for } i<0\}, \\
\B & = & \{X\in\T\ | \ \Hom{\T}{S}{\Sigma^{i}X}=0 \textnormal{ for } i>0\}.
\end{eqnarray*}

Observe that the class of objects $\B$ above coincides with the `torsion class' $\X$ of \eqref{BR_t-structure}. It is natural to ask, therefore, whether analogous theorems to those in \cite{Beligiannis_Reiten} hold in the case of co-$t$-structures. This is indeed the case. In this note, we prove the corresponding results for co-$t$-structures and observe that in the first part of \cite{Co-t-structures} the hypotheses can be relaxed so that $S$ is a connected corigid object (i.e. $\Endo{S}$ is a division ring is not required) and are, in addition, necessary and sufficient.

Recall from \cite{Co-t-structures} that a co-$t$-structure on a triangulated category $\T$ is called {\it non-degenerate} if $\cap_{n\in\Z}\Sigma^{n}\A=\cap_{n\in\Z}\Sigma^{n}\B=0$; following \cite{Beligiannis_Reiten} we shall say it is of {\it finite type} if $\B$ is closed under set indexed coproducts. Consider the following setup.

\setup
{Let $\T$ be a triangulated category with set indexed coproducts. Suppose $\S$ is a set of compact objects in $\T$. Let $\R=\{\Sigma^{i}S\,|\, S\in\S, i<0\}$. Define the following full subcategories of $\T$:
\begin{eqnarray*}
\A & := & \Sigma({}^{\perp}\B); \\
\B & := & \{X\in\T\,|\, \Hom{\T}{S}{\Sigma^{n}X}=0 \textnormal{ for all } S\in\S \textnormal{ and } n>0\}.
\end{eqnarray*}
Note that $\B=\R^{\perp}$ and $\Add{\R}^{\perp}=\R^{\perp}$, where $\Add{\R}$ denotes the smallest full subcategory closed under direct summands of arbitary coproducts of objects of $\R$.}
\label{first_setup}
\esetup

Recall that an additive category $\C$ is called a {\it left triangulated category} if it is equipped with an endofunctor $\Omega:\C\rightarrow\C$ (not necessarily an auto-equivalence) and a class of diagrams $\Omega Z\rightarrow X\rightarrow Y\rightarrow Z$ called {\it left triangle} which satisfy the axioms of a triangulated category except that the left triangles may only be `translated' to the left in (TR2). See \cite[Definitions 2.2 and 2.3]{Beligiannis_Marmardis} for details. {\it Right triangulated categories} are defined similarly.

\lemma
Under the assumptions of Setup \ref{first_setup} we have:
\begin{enumerateroman}
\item $\A$ is a left triangulated subcategory $\T$ which is closed under coproducts and extensions.
\item $\B$ is a right triangulated subcategory of $\T$ which is closed under products and extensions.
\end{enumerateroman}
\label{closure_under_shifts}
\elemma

\pf
This is immediate from Setup \ref{first_setup}; note that in (i) the endofunctor is given by $\Omega=\Sigma^{-1}$ restricted to $\A$, and in (ii) the endofunctor is simply the suspension functor $\Sigma$ restricted to $\B$.
\epf

\

Let $\F$ be a full subcategory of a triangulated category $\T$. A morphism $\phi:X\rightarrow F$ with $F\in\F$ is called a {\it left $\F$-approximation} if the induced morphism $\Hom{\T}{\phi}{F'}:\Hom{\T}{F}{F'}\rightarrow\Hom{\T}{X}{F'}$ is surjective for all $F'\in\F$. Dually, one obtains a {\it right $\F$-approximation}. Left $\F$-approximations are often called {\it $\F$-preenvelopes} and right $\F$-approximations are often called {\it $\F$-precovers}. The full subcategory $F$ is called {\it contravariantly (resp., covariantly) finite} if any object of $\T$ admits a right (resp., left) $\F$-approximation.

\lemma
$\Add{\R}$ is contravariantly finite in $\T$ and for any $X\in\T$ there exists a distinguished triangle $\trilabels{R_{0}}{X}{B_{1}}{f_{0}}{g_{0}}{h_{0}}$ in $\T$ such that
\begin{enumerateroman}
\item $f_{0}$ is a right $\Add{\R}$-approximation of $X$;
\item $0\rightarrow\Hom{\T}{\R}{\Sigma^{n}B_{1}}\rightarrow\Hom{\T}{\R}{\Sigma^{n+1}R_{0}}\rightarrow\Hom{\T}{\R}{\Sigma^{n+1}X}\rightarrow 0$ is a short exact sequence for all $n\geq 0$.
\item The morphism $\Hom{\T}{\Sigma^{-n}g_{0}}{\B}:\Hom{\T}{\Sigma^{-n}B_{1}}{\B}\rightarrow\Hom{\T}{\Sigma^{-n}X}{\B}$ is an isomorphism for all $n>0$ and a surjection for $n=0$.
\end{enumerateroman}
\label{triangle_lemma}
\elemma

\pf
One simply dualises the argument of \cite[Lemma III.2.2]{Beligiannis_Reiten}.
\epf

\thm
Let $\T$ be a triangulated category with set indexed coproducts and $\S$ a set of compact objects of $\T$ as in Setup \ref{first_setup}. Then the pair of full subcategories defined in Setup \ref{first_setup}, 
\begin{eqnarray*}
\A & = & \Sigma({}^{\perp}\B); \\
\B & = & \{X\in\T\,|\, \Hom{\T}{S}{\Sigma^{n}X}=0 \textnormal{ for all } S\in\S \textnormal{ and } n>0\}.
\end{eqnarray*}
defines a co-$t$-structure of finite type on $\T$.
\label{prop1}
\ethm

\pf
We follow the proof of \cite[Theorem III.2.3]{Beligiannis_Reiten}. For an object $X\in\T$, by Lemma \ref{triangle_lemma} we can inductively construct distinguished triangles 
\begin{equation}
\trilabels{R_{n}}{B_{n}}{B_{n+1}}{f_{n}}{g_{n}}{h_{n}}
\label{induction_triangle}
\end{equation}
for $n\geq 0$, where $B_{0}=X$. From these triangles, one obtains a tower of objects and morphisms,
\begin{equation}
X=B_{0}\rightlabel{g_{0}} B_{1}\rightlabel{g_{1}} B_{2}\rightlabel{g_{2}} B_{3}\rightlabel{g_{3}}\cdots\longrightarrow B_{n}\rightlabel{g_{n}} B_{n+1}\longrightarrow\cdots.
\label{tower}
\end{equation}
Recall from \cite{Neeman}, for instance, that the {\it homotopy colimit} of the tower \eqref{tower} is given by the distinguished triangle
$$\trilabel{\coprod_{i=0}^{\infty}B_{i}}{\coprod_{i=0}^{\infty}B_{i}}{\hocolim{B_{i}}}{1-\rm{shift}}{}{}.$$
This induces a morphism $g_{X}:X\to \hocolim B_{n}$.

Note that in the proof of \cite[Theorem III.2.3]{Beligiannis_Reiten} it is shown that this is the reflection of $X$ along the inclusion functor $\iota:\B\rightarrow\T$, and hence one obtains a left adjoint and a $t$-structure. Here, this argument doesn't apply because Lemma \ref{triangle_lemma} is not an exact dual of the corresponding lemma in \cite{Beligiannis_Reiten}. However, it is sufficient to prove that the morphism $g_{X}:X\to \hocolim B_{n}$ is a left $\B$-approximation. First we must verify that $\hocolim{B_{n}}$ is indeed in $\B$.

By construction, the morphism 
$$\Hom{\T}{\R}{\Sigma^{i}(g_{n})}:\Hom{\T}{\R}{\Sigma^{i}B_{n-1}}\rightarrow\Hom{\T}{\R}{\Sigma^{i}B_{n}}$$
is zero for all $n\geq 0$ and all $i>0$. It follows from the short exact sequence
$$0\rightarrow\coprod_{n\geqslant 0}\Hom{\T}{\R}{B_{n}}\rightarrow\coprod_{n\geqslant 0}\Hom{\T}{\R}{B_{n}}\rightarrow\colim\Hom{\T}{\R}{B_{n}}\rightarrow 0$$
that $\colim\Hom{\T}{\R}{B_{n}}=0$. By \cite[Lemma 2.8]{Neeman} we have an isomorphism  $\colim{\Hom{\T}{\R}{B_n}}\cong \Hom{\T}{\R}{\hocolim{B_{n}}}$ since $\R$ consists of compect objects, hence we obtain $\Hom{\T}{\R}{\hocolim B_{n}}=0$ and $\hocolim B_{n}\in\B$.

Now we prove that $g_{X}:X\to \hocolim B_{n}$ is a left $\B$-approximation by following the proof of \cite[Proposition 4.2]{Co-t-structures}. Let $B'\in\B$ and consider the distinguished triangle \eqref{induction_triangle}. By Lemma \ref{triangle_lemma}, the morphism
$$\Hom{\T}{\Sigma^{i}(g_{n})}{B'}:\Hom{\T}{\Sigma^{i}B_{n}}{B'}\rightarrow\Hom{\T}{\Sigma^{i}B_{n-1}}{B'}$$
arising from \eqref{induction_triangle} is an isomorphism for $i<0$ and a surjection for $i=0$. In particular,  given a map $\beta_{0}:X\rightarrow B'$, the fact that we have a surjection for $i=0$ yields the following commutative diagram:
$$\xymatrix{ & B' & & & &  \\
X=B_{0}\ar[r]_-{g_{0}}\ar[ur]^-{\beta_{0}} & B_{1}\ar[r]_-{g_{1}}\ar[u]^-{\beta_{1}} & B_{2}\ar[r]_-{g_{1}}\ar[ul]^-{\beta_{2}} & \cdots\ar[r] & B_{n}\ar[r]_-{g_{n}}\ar[ulll]_-{\beta_{n}} & \cdots \ .}$$
By construction, the composite $$\coprod_{i=0}^{\infty} B_{i}\rightlabel{1-\textnormal{shift}} \coprod_{i=0}^{\infty} B_{i}\rightlabel{\langle\beta_{i}\rangle} B'$$ is zero, where $\langle\beta_{i}\rangle$ is the unique map arising from the coproduct. Thus, we obtain the following commutative diagram
$$\xymatrix{\coprod_{i=0}^{\infty} B_{i}\ar[r]^-{1-\textnormal{shift}}\ar[dr]_-{0} & \coprod_{i=0}^{\infty} B_{i}\ar[r]\ar[d]^-{\langle\beta_{i}\rangle} & \hocolim{B_{i}}\ar[r]\ar@{-->}[dl]^-{\exists} & \Sigma \coprod B_{i} \\
 & B' & & }.$$
It follows that $g_{X}:X\rightarrow\hocolim B_{n}$ is a left $\B$-approximation, as claimed.

Now we need to verify that $(\A,\B)$ as defined in Setup \ref{first_setup} is a co-$t$-structure on $\T$. Conditions $(0)$, $(1)$ and $(2)$ are clear. In order to prove $(3)$ we need the following lemma.

\lemma
The morphism $\Hom{\T}{\Sigma^{i}g_{X}}{\B}:\Hom{\T}{\Sigma^{i}B}{\B}\rightarrow \Hom{\T}{\Sigma^{i}X}{\B}$, where $B=\hocolim B_{n}$, is an isomorphism for $i<0$.
\label{isom}
\elemma

\pflemma
Apply the functor $\Hom{\T}{\Sigma^{i}-}{\B}$ for $i<0$ to the tower \eqref{tower} above to get the inverse tower
\begin{equation}
\cdots\to\Hom{\T}{\Sigma^{i}B_{2}}{\B}\to\Hom{\T}{\Sigma^{i}B_{1}}{\B}\to\Hom{\T}{\Sigma^{i}B_{0}}{\B}.
\label{inverse_tower}
\end{equation}
By Lemma \ref{triangle_lemma}, each morphism in the tower \eqref{inverse_tower} is an isomorphism, so we have $\limit \Hom{\T}{\Sigma^{i}}{\B}\cong\Hom{\T}{\Sigma^{i}B_{0}}{\B}$ for $i<0$. By \cite[Lemma 5.8]{Purity}, there is a short exact sequence
$$\limit^{1}\Hom{\T}{\Sigma^{i+1}B_{n}}{\B}\hookrightarrow\Hom{\T}{\Sigma^{i}\hocolim B_{n}}{\B}\twoheadrightarrow\limit\Hom{\T}{\Sigma^{i}B_{n}}{\B}$$
for $i<0$, where $\limit^{1}$ denotes the first right derived functor of $\limit$; see \cite{Weibel}. For $i\leq 0$, the tower \eqref{inverse_tower} consists of surjective morphisms (for $i<0$, isomorphisms), it particular, it satisfies the Mittag-Leffler condition. It follows by \cite{Weibel} that $\limit^{1}\Hom{\T}{\Sigma^{i+1}B_{n}}{\B}=0$, so that 
$$\Hom{\T}{\Sigma^{i}\hocolim B_{n}}{\B}\rightiso\limit\Hom{\T}{\Sigma^{i}B_{n}}{\B},$$
gives the desired isomorphism.
\epflemma

\

Now we return to the proof of Theorem \ref{prop1} and the verification of condition $(3)$ in definition of a co-$t$-structure. Let $X\in\T$ and take the left $\B$-approximation $g_{X}:X\rightarrow B$, where $B=\hocolim B_{n}$. Extend this morphism to a distinguished triangle:
\begin{equation}
\backtri{A}{X}{B}.
\label{triangle1}
\end{equation}
Applying the functor $\Hom{\T}{-}{\B}$ to \eqref{triangle1} and using Lemma \ref{isom}, one can read off from the resulting long exact sequence that $\Hom{\T}{\Sigma^{-1}A}{\B}=0$, i.e. $A\in\A$, as desired.

Since $\S$ is consists of compact objects, then $\B$ is closed under set indexed coproducts, and hence $(\A,\B)$ is a co-$t$-structure of finite type on $\T$.
\epf

\

In \cite[Definition 4.4.1]{Bondarko} a co-$t$-structure $(\A,\B)$ is called {\it left adjacent} to a (co)-$t$-structure $(\X,\Y)$ if $\B=\X$; cf. the notion of `torsion torsion-free triple' in \cite{Beligiannis_Reiten}. The notion of {\it right adjacency} is defined analogously.

\cor
Under the additional hypotheses 
\begin{enumerate}
\item $\Hom{\T}{S}{\Sigma^{i}S'}=0$ for all $S,S'\in\S$ and for all $i>0$;
\item $\{\Sigma^{i}S\, |\, i\in\Z, S\in\S\}$ is a generating set;
\end{enumerate}
then the co-$t$-structure of Theorem \ref{prop1} is left adjacent to the $t$-structure obtained in \cite[Proposition III.2.8]{Beligiannis_Reiten} (cf. \eqref{BR_t-structure} and Proposition \ref{prop2}).
\label{corollary}
\ecor

\eg
{Consider the homotopy category of spectra $\mathsf{Ho}(\mathsf{Sp})$ and let $\S=\{S^{0}\}$ consist of only the sphere spectrum. Then $S^{0}$ is compact and satisfies the hypotheses of \cite[Proposition III.2.8]{Beligiannis_Reiten} (see \cite{Margolis}, for instance). Thus it follows that the co-$t$-structure induced by Theorem \ref{prop1} is left adjacent to the natural $t$-structure on $\mathsf{Ho}(\mathsf{Sp})$.}
\eeg

\rem
{Theorem \ref{prop1} is established for arbitrary sets of objects $\S$ (whether compact or not) in the case that $\T$ is an `efficient' algebraic triangulated category; see Corollary 3.5 and Proposition 3.15 of \cite{Saorin_Stovicek}.}
\erem

We next obtain an analogue of \cite[Proposition III.2.8]{Beligiannis_Reiten}. We add the following conditions to Setup \ref{first_setup}:

\setup
{Suppose in addition to the conditions satisfied in Setup \ref{first_setup} the set of objects $\S$ also satisfies:
\begin{enumerate}
\item $\Hom{\T}{\Sigma^{i}S}{S'}=0$ for $i>0$ and all objects $S$ and $S'$ in $\S$;
\item $\Hom{\T}{S}{\Sigma S'}=0$ for all objects $S$ and $S'$ in $\S$.
\end{enumerate}}
\label{second_setup}
\esetup

The construction here is slightly different from that of Theorem \ref{prop1} and hinges on the following technical lemma.

\lemma
Under the conditions of Setups \ref{first_setup} and \ref{second_setup}, for any object $X$ in $\T$ there is a left $\B$-approximation $\beta:X\rightarrow B$ with $B\in\B=\R^{\perp}$ such that $\Hom{\T}{S}{\Sigma^{i}(\beta)}:\Hom{\T}{S}{\Sigma^{i}X}\rightarrow\Hom{\T}{S}{\Sigma^{i}B}$ is an isomorphism for $i<1$ and all $S\in\S$.
\label{technical}
\elemma

\pf
Note that in \cite[Proposition 4.1]{Co-t-structures} the hypothesis that $\Endo{S}$ is a division ring (for all $S\in\S$) is not required, and so the lemma follows by \cite[Proposition 4.2 and Lemma 5.1]{Co-t-structures}.
\epf

\prop
Under the conditions of Setups \ref{first_setup} and \ref{second_setup} we have
$$\A \subseteq \{X\in\T \,|\, \Hom{\T}{S}{\Sigma^{i}X}=0 \textnormal{ for } i<0,\, S\in\S\}.$$
Moreover, the following statements are equivalent:
\begin{enumerateroman}
\item The set $\{\Sigma^{i}S\,|\,i\in\Z, S\in\S\}$ is a generating set for $\T$.
\item $\A = \{X\in\T \,|\, \Hom{\T}{S}{\Sigma^{i}X}=0 \textnormal{ for } i<0,\, S\in\S\}$.
\end{enumerateroman}
\label{prop2}
\eprop

\pf
Suppose $(\A,\B)$ is the co-$t$-structure on $\T$ induced in Theorem \ref{prop1}. Let $\bar{\A}= \{X\in\T \,|\, \Hom{\T}{S}{\Sigma^{i}X}=0 \textnormal{ for } i<0,\, S\in\S\}$; in order to show that $\A\subseteq\bar{\A}$, it is sufficient to show that $$\Sigma^{-1}\A\subseteq\Sigma^{-1}\bar{\A}=\{X\in\T \,|\, \Hom{\T}{S}{\Sigma^{i}X}=0 \textnormal{ for } i<1,\, S\in\S\}.$$
Let $X\in\Sigma^{-1}\A$ and consider the left $\B$-approximation $\beta:X\rightarrow B$ arising from Lemma \ref{technical}.  By definition, the induced map $\Hom{\T}{\beta}{B'}:\Hom{\T}{B}{B'}\twoheadrightarrow\Hom{\T}{X}{B'}$ is a surjection for all $B'\in\B$. Since $X\in\Sigma^{-1}\A={}^{\perp}\B$, we have $\Hom{\T}{X}{B'}=0$ for all $B'\in\B$. In particular, setting $B'=B$ yields $\beta=0$. Now by Lemma \ref{technical}, we have the following isomorphism:
$$\Hom{\T}{S}{\Sigma^{i}(\beta)}:\Hom{\T}{S}{\Sigma^{i}X}\rightiso\Hom{\T}{S}{\Sigma^{i}B}$$
for all $S\in\S$ and $i<1$. Since $\beta=0$, it follows that the induced isomorphism $\Hom{\T}{S}{\Sigma^{i}(\beta)}=0$ for all $i<1$, in which case we must have $\Hom{\T}{S}{\Sigma^{i}X}=0$ for all $S\in\S$ and $i<1$, i.e. $X\in\Sigma^{-1}\bar{\A}$, giving the desired inclusion.

$(i)\implies (ii).$ Suppose that $\{\Sigma^{i}S\,|\,i\in\Z, S\in\S\}$ is a generating set for $\T$. We claim that $\A=\bar{\A}$, we have shown the inclusion $\A\subseteq\bar{\A}$ above so we only need to verify that $\bar{\A}\subseteq\A$. Let $X\in\bar{\A}$; by Lemma \ref{technical}, there is a left $\B$-approximation $\beta:\Sigma^{-1}X\rightarrow B$ with $B\in\B$, namely, for any $B'\in\B$ we have a surjection $\Hom{\T}{B}{B'}\twoheadrightarrow\Hom{\T}{\Sigma^{-1}X}{B'}$. We use the argument of \cite[Theorem 5.1]{Co-t-structures}, namely, we have the following isomorphism and equalities:
\begin{eqnarray*}
\Hom{\T}{S}{\Sigma^{i-1}X} & \cong & \Hom{\T}{S}{\Sigma^{i}B} \textnormal{ for all $i<1$ and $S\in\S$ (Lemma \ref{technical})} \\
\Hom{\T}{S}{\Sigma^{i-1}X} & = & 0 \textnormal{ for all $i<1$ and $S\in\S$ (since $X\in\bar{\A}$)} \\
\Hom{\T}{S}{\Sigma^{i}B}   & = & 0 \textnormal{ for all $i>0$ and $S\in\S$ (since $B\in\B$)}.
\end{eqnarray*}
It follows that $\Hom{\T}{S}{\Sigma^{i}B}=0$ for all $i\in\Z$, and thus, since $\{\Sigma^{i}S\,|\,i\in\Z, S\in\S\}$ is a generating set in $\T$, we have $B=0$. Hence $\Hom{\T}{B}{B'}=0$ for all $B'\in\B$, in which case, we have $\Hom{\T}{\Sigma^{-1}X}{B'}=0$. Thus $\Sigma^{-1}X\in{}^{\perp}\B$, i.e. $X\in\A$.

$(ii)\implies (i).$ Let $X$ be an object of $\T$ such that $\Hom{\T}{\Sigma^{i}S}{X}=0$ for all $i\in\Z$ and $S\in\S$. Then $X\in\B$ and $X\in\Sigma^{-1}\A$, i.e. $X\in\Sigma^{-1}\A\cap\B=\{0\}$, whence  $\{\Sigma^{i}S\,|\,i\in\Z, S\in\S\}$ is a generating set for $\T$.
\epf

\rem
{(1) Note that the argument of \cite{Beligiannis_Reiten} does not apply here because one also requires the vanishing of $\Hom{\T}{S}{\Sigma S'}$ for all $S$ and $S'$ in $\S$. To ensure this following the argument of \cite{Beligiannis_Reiten} would require assuming that $\Endo{S}=0$. The condition $\Hom{\T}{S}{\Sigma S'}$ for all $S$ and $S'$ in $\S$ is used in the construction of the left $\B$-approximation in Lemma \ref{technical}.}

\textnormal{(2) Under the hypotheses of Setups \ref{first_setup} and \ref{second_setup}, Proposition \ref{prop2} gives the existence of a co-$t$-structure right adjacent to that obtained in \cite[Theorem III.2.3]{Beligiannis_Reiten}, cf. \eqref{BR_t-structure2}.}
\erem

In the following example we look at the canonical co-$t$-structure induced by a connected cochain differential graded algebra (DGA); for details regarding DGAs see \cite{Bernstein} and \cite{Felix}, also see the motivation presented in \cite{Co-t-structures}. In particular, note that such DGAs arise naturally as the cochain algebras of (simply) connected CW-complexes; see \cite{Felix}.

\eg
{Let $R$ be a connected cochain DGA, i.e $H^{i}(R)=0$ for $i<0$ and $H^{1}(R)=0$, and let $\D(R)$ denote its (unbounded) derived category of differential graded (DG) left $R$-modules. Observe that $H^{i}(M)\cong\Hom{\D(R)}{R}{\Sigma^{i}M}$ for all $i\in\Z$, thus $R$ clearly satisfies the conditions of Proposition \ref{prop2}, and as such one obtains a co-$t$-structure on $\D(R)$.}
\eeg

\

\noindent{\bf Acknowledgement.} The author would like to thank Mar\'ia Jos\'e Souto Salorio for correcting an error in Corollary \ref{corollary}.

\end{document}